\documentclass[reqno]{amsart}
\usepackage{amsmath}
\usepackage{amssymb}
\usepackage{xypic}
\usepackage{mathrsfs}
\usepackage{hyperref}
\usepackage{graphicx}

\date{\today}
          
 \usepackage[usenames,dvipsnames]{pstricks}
\usepackage[shell]{pdftricks}
 \begin{psinputs}
 \usepackage{pst-grad} % For gradients
 \usepackage{pst-plot} % For axes
 \usepackage{epsfig}
 \end{psinputs}

\begin{document}
\title{Transversality Theorems for the weak topology}
\author{Saurabh Trivedi} 
\address{LATP (UMR 6632), Centre de Math\'ematiques et Informatique, Universit\'e de provence, 39 rue Joliot-Curie, 13013 Marseille, France.}

\email{saurabh@cmi.univ-mrs.fr}
\maketitle

\parskip .20cm

\newtheorem{thm}{Theorem}[section]
\newtheorem{lem}[thm]{Lemma}

\theoremstyle{theorem}
\newtheorem{prop}[thm]{Proposition}
\newtheorem{definition}[thm]{Definition}
\newtheorem{example}[thm]{Example}
\newtheorem{xca}[thm]{Exercise}

\theoremstyle{remark}
\newtheorem{rem}[thm]{Remark}

\numberwithin{equation}{section}

\newenvironment{pf}{\noindent \textbf{Proof:}}{\hfill $\square$\\}
\newcommand{\bb}{\mathbb}
\newcommand{\al}{\mathcal}
\newcommand{\ak}{\mathfrak}
\newcommand{\fs}{\mathscr}

\renewcommand{\thefootnote}{\fnsymbol{footnote}}

\begin{abstract}
In his paper \cite{Trotman}, Trotman proves, using the techniques of the Thom transversality theorem, that under some conditions on the dimensions of the manifolds under consideration, openness of the set of maps transverse to a stratification in the strong (Whitney) topology implies that the stratification is $(a)$-regular.  Here we first discuss the Thom transversality theorem for the weak topology and then give a similiar kind of result for the weak topology, under very weak hypotheses.  Recently several transversality theorems have been proved for complex manifolds and holomorphic maps (see \cite{Kaliman} and \cite{Forstneric}). In view of these transversality theorems we also prove a result analogous to Trotman's result in the complex case.
\end{abstract}

\section{Smooth case}
Let $M$ and $N$ be smooth manifolds. Denote by $C^{\infty}(M,N)$, the set of all smooth maps between $M$ and $N$. The set $C^{\infty}(M,N)$ can then be given two topologies, the weak topology and the strong (Whitney) topology (see page 35 in Hirsch \cite{Hirsch}  for the definitions, see also \cite{Plessis} for more about the function space topologies). Many authors prefer to use a definition of this topology on the function spaces via jets; see page 42 in Golubitsky and Guillemin \cite{Golubitsky} and also a somewhat detailed discussion about the different topologies on function spaces in du Plessis and Vosegaard \cite{Plessis1}. We will follow the approach of Hirsch \cite{Hirsch}.  Denote by $C^{\infty}_W(M,N)$ and $C^{\infty}_S(M,N)$, the space of smooth maps between $M$ and $N$ with the weak topology and the strong topology respectively. We prove the results for the weak topology and state them for the strong topology.

We say that a smooth map $f:M \rightarrow N$ is transverse to a submanifold $S \subset N$ at $x \in M$, denoted $f \pitchfork_x S$, if either $f(x) \notin S$ or $f(x) \in S$ and $T_{f(x)}S + Df_x(T_xM) = T_{f(x)}N$. A map $f$ is transverse to a submanifold on a subset $K$ of the domain, denoted by $f  \pitchfork_K \Sigma$, if it is transverse at all points of the subset $K$. Notice that if the codimension of $S$ is greater than the dimension of $M$ then a map $f : M \rightarrow N$ is transverse to $S$ if and only if $f(M) \cap S = \emptyset$, i.e., if the image of $M$ under $f$ is disjoint from $S$. 

A major result that describes the transversal intersection property of smooth maps with respect to the submanifolds of the target manifold is the Thom transversality theorem (theorem 2.1 on page 74 in \cite{Hirsch}). More precisely,
\pagebreak
\begin{thm}\label{trans1} Let $M$ and $N$ be smooth manifolds, $S \subset N$ a submanifold. Then,

(a) $T_S=\{f \in C^{\infty}(M,N) : f \pitchfork S\}$ is a dense subset of $C^{\infty}_W(M,N)$ as well as of $C^{\infty}_S(M,N)$.

(b) Suppose $S$ is closed in $N$ and $K \subset M$. Then $\{f \in C^{\infty}(M,N) : f \pitchfork_K S\}$ is open in $C^{\infty}_W(M,N)$ if $K$ is compact and open in $C^{\infty}_S(M,N)$ if $K$ is closed.
\end{thm}

What would seem to be an obvious generalization of the transversality theorem \ref{trans1} is to replace the submanifold $S$ by a collection of submanifolds such that their union is a closed subset, which appears as exercise 3 on page 59 in Golubitsky and Guillemin's book \cite{Golubitsky}. Unfortunately the exercise is not correct as stated and we will provide a counterexample later. A similar kind of mistake appears in exercise 8 on page 83 in Hirsch's book \cite{Hirsch}; we will also give a counterexample to this exercise.  

In fact there seem to be no complete correct statements published to date of such results and perhaps due to the mistakes in the standard textbooks, many recent papers contain errors of a similar kind. See Trotman \cite{Trotman} for a brief discussion about some mistakes occurring in papers of Thom, Chenciner and Wall. More recently, Loi \cite{Loi} proves  a transversality theorem for definable maps in $o$-minimal structures, but his definition of the definable jet bundle is incorrect and due to the problems with definitions his proof is not correct as stated. Also, it is generally believed that the openness of the set of maps transverse to a stratification in the strong topology implies $a$-regularity of the stratification without any restrictions on the dimensions of manifolds, see for example remark 1.3.3 on page 38 in Goresky and Macpherson \cite{Goresky}, but this is not true; see remark \ref{rem1} below. In fact Goresky and Macpherson \cite{Goresky} also prove a more general transversality theorem, proposition 1.3.2 on page 38 and claim that the openness of maps, restricted to a stratification in the source manifold, transverse to a stratification in the target manifold in the strong topology implies that both stratifications are $a$-regular and refer to Trotman \cite{Trotman}, but Trotman does not prove such a statement in his paper. 

Here we prove a proposition which is more general than the transversality theorem and also gives the correct formulation of the incorrect exercises mentioned above. We need the following definitions: 

A {\bf stratification} $\Sigma$ of a subset $V$ of a manifold $M$ is a locally finite  partition  of $V$ into submanifolds of $M$. The submanifolds in the partition are called strata. By a locally finite partition we mean that each point of $V$ has a neighbourhood meeting only finitely many strata.

Let $S_1$ and $S_2$ be two strata of $\Sigma$, $S_2$ is said to be {\bf ${\textbf a}$-regular} over $S_1$ at $x \in S_1 \cap \overline{S_2}$ if for every sequence of points $\{y_i\}$ in $S_2$ converging to $x$ such that $\lim_{i\rightarrow \infty} T_{y_i} S_2$ exists, we have 
$$\lim_{i\rightarrow \infty} T_{y_i} S_2= \tau \,\,\,\Rightarrow \,\,\,T_xS_1 \subset \tau.$$ 
A stratification is called $a$-regular if every pair of strata $(S_i,S_j)$ is $a$-regular at every point in the intersections $S_i \cap \overline{S_j}$ and $S_j \cap \overline{S_i}$.

\begin{figure}[h]
\includegraphics[scale=.85]{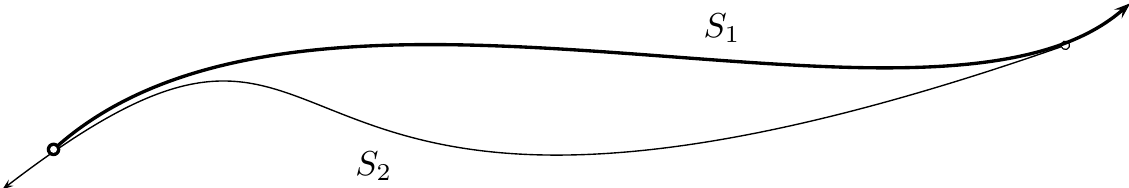}
\caption{}\label{fig1}
\end{figure}

Notice that our definition of a stratification allows that none of the strata be a closed subset; see figure \ref{fig1}. Both $S_1$ and $S_2$ are not closed and yet their union is closed.

We say that a map $f$ is transverse to a stratification $\Sigma$ at a point $x$ if $f$ is transverse to each stratum $S_i$ at $x$.
 
Now we state the result:

\begin{prop} Let $M$ and $N$ be smooth manifolds and let $\Sigma$ be an $a$-regular stratification of a closed subset $S \subset N$. Then,

(a) $T_{\Sigma}=\{f \in C^{\infty}(M,N) : f \pitchfork \Sigma\}$ is a dense subset of both $C^{\infty}_W(M,N)$ and $C^{\infty}_S(M,N)$. 

(b) Suppose $K \subset M$. Then $T_K = \{f \in C^{\infty}(M,N) : f \pitchfork_K \Sigma\}$ is open in $C^{\infty}_W(M,N)$ if $K$ is compact and open in $C^{\infty}_S(M,N)$ if $K$ is closed.
\end{prop}

\begin{rem}
The strata of a stratification need not be closed submanifolds and they can be of same dimension but the union of the strata should be closed.

The denseness of transversal maps is proved in Hirsch \cite{Hirsch} (theorem 2.5 on page 78 and theorem 2.8 on page 80). In fact, we don't need $a$-regularity to prove the denseness. Here we prove results about openness of transversal maps.

The transversality theorem can be seen as a special case of our proposition, namely, when the stratification of the underlying closed set has only one stratum. 
\end{rem}

\begin{lem}\label{lem1}
Let $M$ and $N$ be manifolds and let $\Sigma$ be an $a$-regular stratification of a closed subset $S$ of $N$ and let $f \in C^{\infty}(M,N)$. If $x \in M$ is such that $f \pitchfork_x \Sigma$ then there exists a coordinate chart $(\phi, U)$ of $M$ at $x$ such that for each compact $K \subset U$ there is a weak subbasic  neighbourhood $\mathcal N(f,(\phi, U),(\psi, V),K,\epsilon)$\footnote{Set of all maps $g \in C^{\infty}(M,N)$ such that $g(K) \subset V$ and $||D^rf_{\phi,\psi}(x) - D^rg_{\phi,\psi}(x)|| < \epsilon$ for all $r \geq 0$ and all $x \in \phi(K)$} of $f$ such that every member $g$ of this neighbourhood satisfies $g \pitchfork_K \Sigma$.
\end{lem}

\begin{pf}
Without loss of generality we can assume that $\Sigma$ has only two distinct strata $S_1, S_2$. 
We have four cases:

{\bf 1.}  $\boldsymbol {f(x) \notin S}$:

Since $S$ is closed, we can choose a coordinate chart $(\psi, V)$ around $f(x)$ such that $V \cap S = \emptyset$. But since $f^{-1}(V)$ is open containing $x$, we can choose a chart $(\phi, U)$ around $x$ such that $f(U) \subset V$.

Thus for each compact set $K \subset U$ we can choose $\epsilon >0$ small enough  and take a weak subbasic neighbourhood   $\mathcal N(f,(\phi, U),(\psi, V),K,\epsilon)$ of $f$ such that if $g$ belongs to this neighbourhood then in particular $g(K) \subset V$; hence $g(K) \cap S = \emptyset$ and thus $g \pitchfork_K \Sigma$.

{\bf 2. $\boldsymbol{f(x) \in S}$ but $\boldsymbol{f(x) \notin S_1 \cap \overline{S_2}}$ and $\boldsymbol{f(x) \notin S_2 \cap \overline{S_1}}$:}

Suppose here that the dimensions of $M$ and $N$ are $m$ and $n$ respectively. Let the dimension of $S_1$ be $s$.

Since $f(x) \in S$, there exists a stratum, say $S_1$, with $f(x) \in S_1$. Let $(\psi,V)$ be a chart for $(N,S_1)$ at $f(x)$ such that $V \cap S_2 = \emptyset$, this is possible since $f(x) \notin S_1 \cap \overline{S_2}$ and $f(x) \notin S_2 \cap \overline{S_1}$. Let $(\phi,U')$ be a chart for $M$ at $x$ such that $f(U') \subset V$.  Let $f_{\phi, \psi} =  \psi \circ f \circ \phi^{-1}$.  Let $\pi : \bb R^n \rightarrow \bb R^n/(\bb R^s \times 0)$ be the natural projection map. The situation is described below:
\[
\xymatrix{
U' \ar[r]^{\phi} \ar[d]_f & \bb R^m \ar[d]^{f_{\phi,\psi}} \ar[rrd]&\\
 V \ar[r]^{\psi}& \bb R^n \ar[rr]^{\pi~~~~~~~~~~~~} & &\bb R^n/ \bb R^s \times 0
}
\]

Denote by $L(\bb R^m,\bb R^n/(\bb R^s \times 0))$ the set of linear maps between $\bb R^m$ and $\bb R^n/(\bb R^s \times 0)$; it is a normed space. 

Define a map $\eta : U' \rightarrow \bb R$ such that $\eta(u)$ denotes the distance of the linear map 
$$D_{\phi(u)} ( \pi \circ f_{\phi,\psi}) : \bb R^m \rightarrow (\bb R^n/\bb R^s \times 0)$$
from the set of nonsurjective linear maps in $L(\bb R^m,\bb R^n/(\bb R^s \times 0))$.

Note that $\eta$ is a continuous map and $\eta(x) >0$ since $f$ is transverse to $S_1$ at $x$ and the nonsurjective linear maps $L(\bb R^m,\bb R^n/(\bb R^s \times 0))$ form a closed set. Thus, there exists a neighbourhood $U$ of $x$ on which $\eta >0$. Hence if $K$ is any compact subset of $U$ then $\eta$ will assume a positive minimum value, say $2 \epsilon$, on $K$.  Clearly then, the weak neighbourhood $\mathcal N(f,(\phi|_U, U),(\psi, V),K,\epsilon)$ is the required weak neighbourhood of $f$. 

{\bf 3. $\boldsymbol {f(x) \in S}$ and $\boldsymbol{f(x) \in S_1 \cap \overline{S_2}}$:}

Using the same argument as in the case 2, we can find coordinate charts $(\phi'',U'')$ around $x$ and $(\psi'',V'')$ around $f(x)$ such that for each compact $K \subset U''$ there is a weak subbasic neighbourhood $\mathcal N(f,(\phi'', U''),(\psi'', V''),K,\epsilon'')$ such that every member of this neighbourhood is transverse to $S_1$ on $K$.

Now, suppose that for each neighbourhood $U$ of $x$, there is a compact set $K \subset U$ such that every neighbourhood $\mathcal N(f,(\phi, U),(\psi, V),K,\xi)$ contains a map $g$ satisfying $g \not\pitchfork_K S_2$.

Choose $\{U_i\}_{i=1}^{\infty}$  to be a basis for the neighbourhoods of $x$. Then, for each $i$ there is a compact set $K_i \subset U_i$, a point $y_i \in K_i$ and a map $g_i \in \mathcal N(f,(\phi|_{U_i}, U_i),(\psi, V)$ $,K_i,1/i)$ such that $g_i \not\pitchfork_{y_i} S_2$, where $(\phi,U),(\psi,V)$ are fixed charts for $M$ and $N$ at $x$ and $f(x)$ respectively.

Note first that for each $i$ 
$$|\psi g_i(y_i) - \psi f(y_i)| < 1/i$$
and also that
$$|\psi f(y_i) - \psi f(x)| < \epsilon_i$$
 where $\epsilon_i \rightarrow 0$ as $i \rightarrow \infty$. Hence, by the triangle inequality, $g(y_i) \rightarrow f(x)$ as $i \rightarrow \infty$.
 
 Since $g_i \not\pitchfork_{y_i} S_2$ it follows that 
\begin{eqnarray*} 
& \dim T_{g_i(y_i)} N & > \dim\left (T_{g_i(y_i)}S_2 + D_{y_i}g_{y_i}(T_{y_i}M)\right )\\
\Rightarrow & \dim T_{f(x)}N & >  \dim \left (T_{g_i(y_i)}S_2 + D_{y_i}g_{y_i}(T_{y_i}M)\right )
\end{eqnarray*}

Taking limits\footnote{We may assume that all of the above limits exist by taking subsequences if necessary.} on both sides and using the properties of sequences of points in grassmannians we have:
\begin{eqnarray*} 
\dim T_{f(x)}N & > & \lim_{i\rightarrow  \infty} \dim \left (T_{g_i(y_i)}S_2 + D_{y_i}g_{y_i}(T_{y_i}M)\right )\\
& = & \dim \lim_{i\rightarrow \infty} \left (T_{g_i(y_i)}S_2 + D_{y_i}g_{y_i}(T_{y_i}M)\right )\\
& \geq & \dim \left (\lim_{i\rightarrow \infty} T_{g_i(y_i)}S_2 + \lim_{n\rightarrow \infty}  D_{y_i}g_{y_i}(T_{y_i}M)\right )\\
 &\geq & \dim \left (\lim_{i\rightarrow \infty} T_{g_i(y_i)}S_2 +  D_xf(T_xM)\right )
 \end{eqnarray*}

But since $S_2$ is $(a)-$regular over $S_1$ at $f(x)$ we have
$$\lim_{i\rightarrow \infty} T_{g_i(y_i)} S_2 \supset T_{f(x)} S_1$$

Thus it follows that
$$\dim T_{f(x)}N >  \dim\left ( T_{f(x)}S_1 + D_xf(T_xM)\right)$$
which is a contradiction to the fact that $f \pitchfork_x \Sigma$. Thus, there exists a chart      $(\phi',U')$ around $x$ and $(\psi',V')$ around $f(x)$ such that for each compact $K \subset U'$ the subbasic neighbourhood of $f$, $\mathcal N(f,(\phi', U'),(\psi', V),K,\epsilon')$ has the property that all its members are transverse to $S_2$ on all of $K$. 

Set $U = U'\cap U''$, $V = V'$, $\phi = \phi'|_U$ and $\psi = \psi'$. It is easy to see that for a suitable $\epsilon$ and any compact $K \subset U$, the subbasic neighbourhood $\mathcal N(f) = \mathcal N(f,(\phi, U),(\psi, V),K,\epsilon)$ satisfies,
$$\mathcal N(f) \subset \mathcal N(f,(\phi', U'),(\psi', V'),K,\epsilon') \cap\mathcal N(f,(\phi'', U''),(\psi'', V''),K,\epsilon'')$$
and all its members are transverse to $S_1$ and $S_2$ on $K$.

{\bf 4. $\boldsymbol{f(x) \in S}$ and $\boldsymbol{f(x) \in S_2 \cap \overline{S_1}}$:}

The proof is exactly the same as the case 3 with $S_1$ and $S_2$ interchanged.
\end{pf}

\noindent
{\bf Proof of proposition 1.2.}
Let $f \in T_K$. To prove that $T_k$ is open, we show that there exists a weak open neighbourhood of $f$ which is contained in $T_K$. Since $f$ is transverse to $\Sigma$ at each $x \in K$, by lemma \ref{lem1} , for each $x \in K$ there exists a chart $U_x$ with the property that for each compact set $K_x \subset U_x$ there is a neighbourhood 
$N(f,(\phi_x, U_x),(\psi_x, V_x),K_x,\epsilon_x)$ such that each member of this neighbourhood is transverse to $\Sigma$ on all of $K_x$.  Since $K$ is compact, we can choose a finite subcollection $\{U_{x_1},\ldots,U_{x_r}\}$ of the coordinate neighbourhoods $\{U_x\}_{x \in L}$,  such that $K \subset \cup_{i=1}^r K_{x_i}$. But then the intersection 
$$\cap_{i=1}^r N(f,(\phi_{x_i}, U_{x_i}),(\psi_{x_i}, V_{x_i}),K_{x_i},\epsilon)$$ 
($\epsilon = \min\{\epsilon_{x_i}\}$)
is a weak open neighbourhood of $f$ and is contained in $T_K$, as required.

The corresponding statement about the strong topology can be proved by taking a countable covering of the closed set $K$ and then taking the countable intersection of the weak neighborhoods thus obtained. \hfill{$\square$}

\begin{rem}\label{rem1}
If, in a stratification, the codimension of a stratum of the minimal dimension is greater than the dimension of $M$, then we don't need $a$-regularity to prove that the set of maps transversal to the stratification on a compact set is open.

Our result is more general than the result of exercise 15 on page 84 in \cite{Hirsch}, for we don't need any of the submanifolds of the stratification to be closed and they can be of same dimension.
\end{rem}

Consider the following examples:

{\bf 1.} Exercise 8(a) on page 83 in Hirsch \cite{Hirsch} amounts to the following:

\emph{Let $M$ and $N$ be smooth manifolds and let $A_0,\dots, A_q$ be submanifolds of $N$ such that $\cup A_i$ be compact and the $A_i's$ form a submanifold complex , then the set of maps transverse to each of $A_i$ is open in the weak topology.} 

\begin{figure}[h]
\includegraphics[scale=.85]{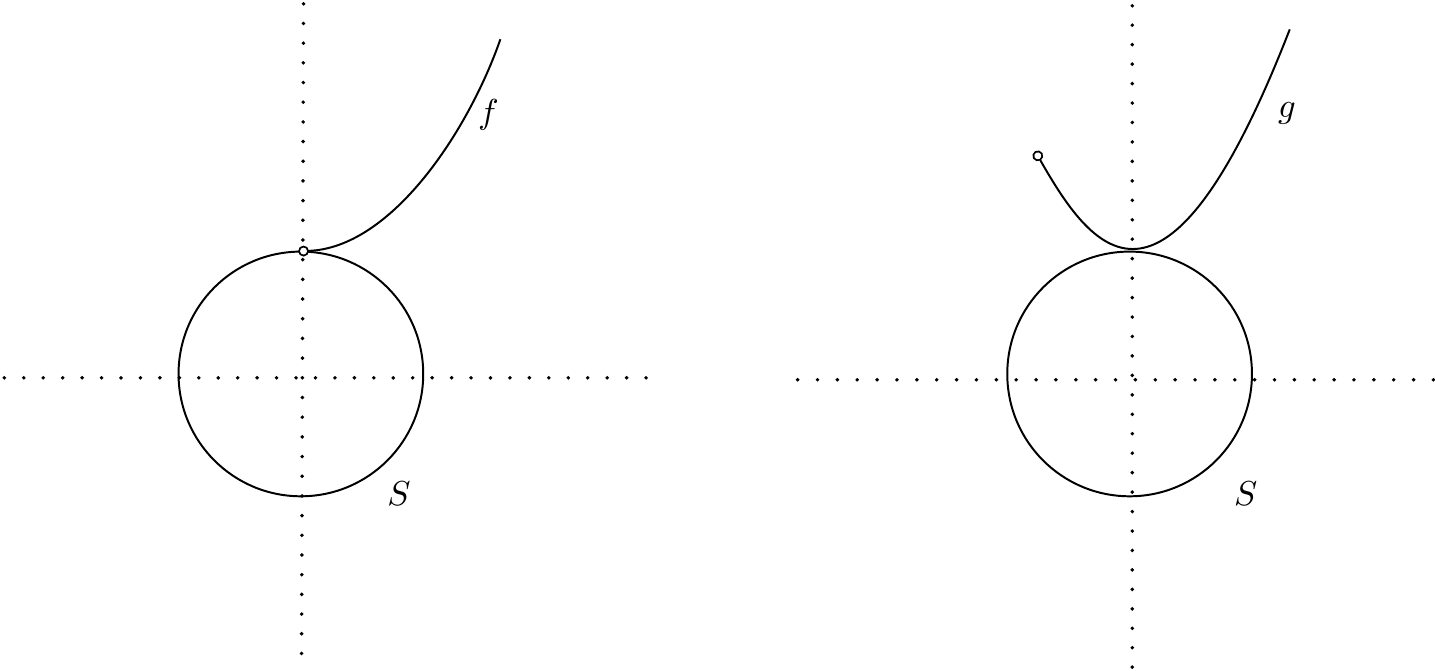}
\caption{}\label{fig2}
\end{figure}

Now, let $M = (0,\infty)$, $N = \bb R^2$ and $S = \{(x,y) \in N  : x^2 + y^2 = 1\}$.  Let $f : M \rightarrow N$ be defined by $f(x) = (x, x^2+1)$. Then $f$ is transverse to $S$ at all points of $M$ but by an arbitrary perturbation of $f$ of the type $g(x) = (x-c, (x-c)^2 + 1)$ for $c>0$, we can find a map which is not transverse to $S$; see figure \ref{fig2}. 

In the above example $S$ is a compact submanifold of $N$, but there exist maps arbitrarily close to $f$ in the weak topology, which are not transverse to $S$.  So, the set of maps transverse to $S$ in the above case is not open in the weak topology. Thus, it is a counter example to the exercise 8 on page 83 in \cite{Hirsch}.

{\bf 2.} Let $M = \bb R$, $N = \bb R^2$ and $S = S_1 \cup S_2 = \bb R^+ \times 0 \cup 0 \times \bb R$. Then $S$ is a closed subset of $N$ and, $S_1$ and $S_2$ are submanifolds. Let $f: M \rightarrow N$ be defined by $f(x) = (x,x^2)$. 

\begin{figure}[h]
\includegraphics[scale=.85]{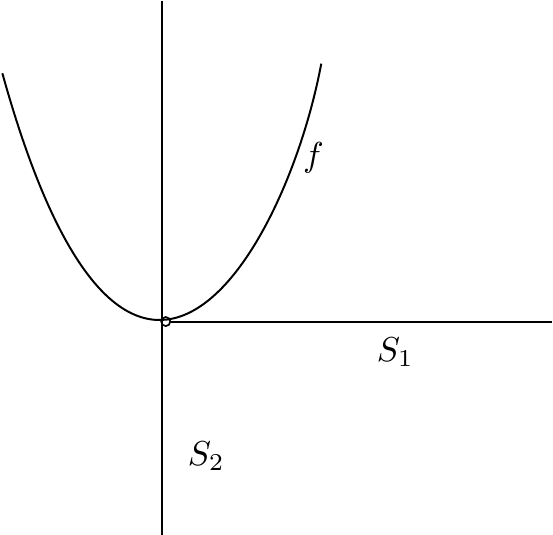}
\caption{}\label{fig3}
\end{figure}

Notice that $f$ is transverse to $S_1$ as well as $S_2$ but we can find maps arbitrarily close to $f$ not transverse to $S_1$ by simply shifting $f$ to the right side; see figure \ref{fig3}. Thus this is a counterexample to exercise 3 on page 59 in \cite{Golubitsky}.\\

Now, let $M$ and $N$ be manifolds and $\Sigma$ be a stratification of a closed subset $V$ of $N$. Suppose that the codimension of a stratum of minimal dimension is less than or equal to the dimension of $M$. Then, in fact, we can prove the converse of the statement of our proposition. A similar kind of result is proved in \cite{Trotman} for the strong topology. 

We suppose that the strata are of dimension at least $1$.

\begin{prop} \label{main} Let $M$ and $N$ be smooth manifolds and $\Sigma$ be a stratification of a closed subset $V$ of $N$. Let $r = \min \{\dim S_i : S_i \,\, \text{is a stratum in} \,\, \Sigma\}$, $\dim N = n$ and $ \dim M \geq n - r$.  If for some $ w \in M$, the set $T_{\Sigma} = \{f \in C^{\infty}(M,N) : f \pitchfork_w \Sigma\}$ is open in the weak topology, then $\Sigma$ is an $a$-regular stratification. 
\end{prop}

\begin{pf} Let $X$ and $Y$ be two distinct strata in $\Sigma$ such that $Y$ is not $a$-regular over $X$ at some $x \in X \cap \overline{Y}$. Then, there exists a sequence $\{y_i\}_{i=1}^{\infty}$ in $Y$ such that $\lim_{i \rightarrow \infty} y_i = x$, $\lim_{i \rightarrow \infty} T_{y_i} Y = \tau$ but $T_xX \not\subset \tau$.

Let $v \in T_x X$ be such that $v \notin \tau$ and $E$ be the one dimensional subspace of $T_x N$ spanned by $v$. Now, choose a basis for $T_xN$ such that we can write
\begin{eqnarray*}
T_x X & =&  E \oplus W_1 \oplus T_1\\
\tau & = & T_1 \oplus T_2\\
T_x N &=& E \oplus W_1 \oplus W_2 \oplus T_1 \oplus T_2
\end{eqnarray*}
where $T_1$, $T_2,$ $W_1,$ $W_2$ are subspaces of $T_xN$, $T_1 = T_x X \cap \tau$. Then find a subspace $H$ of $T_x N$ with $\dim H = \dim N - r$, such that
$$T_2 \oplus W_2 \subseteq H \subseteq T_1 \oplus T_2 \oplus W_1 \oplus W_2.$$
Then, we have
\begin{equation}
H + T_x X = T_x N 
 \end{equation}
and
\begin{equation}
 H + \tau \neq T_x N.
\end{equation}
Let $(\phi, U)$ be a coordinate chart around $w \in M$ such that $\phi(w) = \overrightarrow{0} \in \bb R^{m}$ and $(\psi, V)$ be a coordinate chart around $x \in N$ such that $\psi(x) = \overrightarrow{0} \in \bb R^n$.  

Then, $D_{x}\psi : T_x N \rightarrow \bb R^n$ is a linear isomorphism of vector spaces.  Under this isomorphism (1.1) and (1.2) become
\begin{equation}
D_x\psi(H) + D_x\psi(T_x X) = \bb R^n 
 \end{equation}
and
\begin{equation}
 D_x\psi(H) + D_x\psi(\tau) \neq \bb R^n.
\end{equation}

Let $D_x \psi(v) = {v'}$. Notice that ${v'}  \not\in D_x\psi(H) + D_x \psi(\tau)$. Suppose that the dimension of $H + \tau$ is $p$. Then $n-r \leq p < n$. Let $\{{v_1},\dots,{ v_{p}}\}$ be a basis of $D_x\psi(H) + D_x\psi(\tau)$ such that $\{v_1,\ldots,v_{n-r}\}$ forms a basis of $H$. Extend this basis to a basis $\{{v_1},\ldots,{v_{p}},{v_{p+1}},\ldots,{v_{n-1}},{v'}\}$ of $\bb R^n$.

Let $m$ be the dimension of $M$ and define a map $L : \bb R^{m} \rightarrow \bb R^n$ (this map is well defined because $m \geq n-r$) by, 
$$L(a_1,\ldots,a_m) = a_1 {v_1} + a_2 {v_2} + \ldots + a_{n-r} {v_{n-r}}.$$
Notice that, $L(\bb R^{m}) = D_x\psi(H)$.

By using a bump function we may construct a map $f \in C^{\infty}(M,N)$ such that $f_{\phi,\psi} = \psi \circ f \circ \phi^{-1} = L$ on some neighbourhood of ${0} \in \bb R^m$ and thus $f(w)=x$. Since $L$ is a linear map, $D_{\phi(w)}f_{\phi,\psi} = L$ on the same neighbourhood. 

The situation is described below:
\[
\xymatrix{
U \ar[r]^{\phi} \ar[d]_f & \bb R^m \ar[d]^{f_{\phi,\psi} = L} \\
 V \ar[r]^{\psi}& \bb R^n }
\]

Thus by (1.3),
\begin{eqnarray*}
D_{\phi(w)} f_{\phi,\psi}(D_{w}\phi(T_wM))  + D_x{\psi}(T_xX) & = & D_{\phi(w)} f_{\phi,\psi}(\bb R^m)  + D_x{\psi}(T_xX) \\
& = & L(\bb R^m) + D_x{\psi}(T_xX)\\
& = & D_x\psi(H) + D_x\psi(T_xX)\\ 
&=& \bb R^n.
\end{eqnarray*}

This implies that, $f \pitchfork_w X$ and $f \pitchfork_w Y$ (since $x \notin Y$).

We will now construct a sequence of smooth maps $\{f^k\}$ which converge to $f$ in the weak topology such that for sufficiently large $k$, $f_k \not\pitchfork_w \Sigma$, which will be a contradiction to the hypothesis that the set $T_{\Sigma} = \{f \in C^{\infty}(M,N) : f \pitchfork_w \Sigma\}$ is open in the weak topology.

First note that $y_k \in V$ for sufficiently large $k$. Thus 
$$\lim_{k \rightarrow \infty} D_{y_k}\psi(T_{y_k} Y) = D_x\psi(\tau).$$
 Now, choose a basis $\{v_1^k,v_2^k,\ldots,v_{p}^k,$ $v_{p+1}^k,\ldots,$ $v_{n-1}^k, v^k\}$ of $D_{y_k}\psi(T_{y_k}N)$ such that $D_{y_k}(T_{y_k}Y)$ belongs to the span of $\{v_1^k,\ldots,v_{p}^k\}$ and 
\begin{equation}
\lim_{k \rightarrow \infty} v_i^k = {v_i} \,\,\text{and}\,\, \lim_{k \rightarrow \infty} v^k = {v'}.
\end{equation}
Let $H^k$ be the subspace of $D_{y_k}\psi(T_{y_k}N)$ spanned by $v_1^k,\ldots,v_{n-r}^k$. Then, by (1.5) we have $\lim_{k\rightarrow \infty} H^k = D_x{\psi}(H)$ and 
\begin{equation}
H^k + D_{y_k}\psi(T_{y_k}Y) \neq \bb R^n
\end{equation}
since the left hand side of (1.6) is spanned by a subset of $\{v_1^k,\ldots,v_p^k\}$ and $p < n$.

Now, define the map $f^k \in C^{\infty}(M,N)$ first on the open set $U \subset M$ by the following formula for its local representative:
$$f_{\phi,\psi}^k(z) = f_{\phi,\psi}(z) + \lambda(z) \left ( \psi(y_k) +  z_1({v_1^k} -{v_1}) + \ldots +z_{n-r} ({v_{n-r}^k} - {v_{n-r}})\right ).$$
where $\lambda : \phi(U) \rightarrow \bb R$ is a smooth bump function defined by
\begin{equation*}
\lambda = 
\begin{cases} 1& \text{on a compact neighbourhood $K$ of $\overline{0}$}\\
0 & \text{outside  a relatively compact neighbourhood of $K$ contained in $\phi(U)$}
\end{cases}
\end{equation*}

Notice that $f^k_{\phi,\psi}(w) = y_k$ and $D_{\phi(w)}f_{\phi,\psi}^k$ is the linear map $L^k : \bb R^m \rightarrow \bb R^n$ defined by
$$L^k(z) = z_1 v_1^k + \ldots +z_{n-r} v_{n-r}^k.$$
and so by (1.6),
\begin{eqnarray*}
D_{\phi(w)} f^k_{\phi,\psi}(D_{w}\phi(T_wM))  + D_{y_k}{\psi}(T_{y_k}Y) & = & D_{\phi(w)} f^k_{\phi,\psi}(\bb R^m)  + D_{y_k}{\psi}(T_{y_k}Y) \\
& = & L^k(\bb R^m) + D_{y_k}{\psi}(T_{y_k}Y)\\
& = & H^k + D_{y_k}\psi(T_{y_k}Y)\\ 
&\neq& \bb R^n.
\end{eqnarray*}

Now define $f^k$ outside $U$ to be equal to $f$. This gives a smooth map $f^k$. Notice that for sufficiently large $k$, by (1.6), $f^k \not\pitchfork_w Y$, which is the same as saying, $f^k \not\pitchfork_w \Sigma$ . It is easy to see that $\lim_{k\rightarrow \infty} f^k = f$ in the weak topology (in fact $f^k$ tends to $f$ also in the strong topology), which is a contradiction to the hypothesis that $T_{\Sigma} = \{f \in C^{\infty}(M,N) : f \pitchfork_w \Sigma\}$ is open in the weak topology.
\end{pf}

\begin{rem}
Our proof is inspired by the proof of the main theorem in \cite{Trotman} for the strong topology. We have explicitly constructed the maps to arrive at the contradiction while in \cite{Trotman} Trotman only mentions the existence of such maps.
\end{rem}

\section{Holomorphic case}

Our next aim is to prove, if possible, a similar result for complex manifolds and holomorphic maps between them.

Let $M$ and $N$ be complex manifolds and let $\al H(M,N)$ be the set of all holomorphic maps between $M$ and $N$. 

If $M$ is compact then any holomorphic map between $M$ and $N = \bb C^n$ (for any positive integer $n$) must be constant. Let $x \in N$ and $f_x : M\rightarrow N$ be the constant map taking all elements of $M$ to $x$. Then, we have the bijection $f_x \mapsto x$ between $\al H(M,N)$ and $N$.  

Now, regard $M$ and $N$ as $C^1$ manifolds and denote by $C^1(M,N)$ the set of all $C^1$ maps between $M$ and $N$. Give $C^1(M,N)$ the weak topology (the weak subbasic open neighbourhoods of the weak topology contain functions which are close together with their derivatives on some compact subset of $M$ (see page 35 in \cite{Hirsch} for definitions)). 

Choose a metric on $N$ and give it the metric topology. The relative topology on $\al H(M,N)$ under the weak topology of $C^1(M,N)$ is same as the smallest topology on $\al H(M,N)$ which makes the above bijection a continuous map. We call this topology the weak topology on $\al H(M,N)$.

Now, let $S \subset N$ be a complex submanifold  of $N$. Then, a holomorphic map $g_x : M \rightarrow N$ ($M$ compact) is transverse to $S$, denoted $f \pitchfork_m S$ at any point $m \in M$ if and only if $x \notin S$.  Thus we have the following trivial result:

\begin{prop} Let $M$ be a compact complex manifold. Let $S \subset \bb C^n$ be a complex submanifold. Then, the set $T_S = \{f \in \al H(M,\bb C^n) : f \pitchfork S\}$ is open in $\al H(M,\bb C^n)$ with the weak topology if $S$ is closed.
\end{prop}

We can extend the above result in the following way.

Let $X$ be a complex analytic subvariety of $N$. Then, by a result of Whitney \cite{Whitney}, $X$ can be stratified into complex submanifolds. In fact, $X$ can be stratified such that the strata fit together well in the sense that their tangent spaces satisfy nice regularity conditions. But we don't need these regularity conditions in the following. We have,

\begin{prop} Let $M$ be a compact complex manifold. Let $X \subset \bb C^n$ be a complex analytic subvariety. Let $\Sigma$ be a stratification of $X$. Then the set  $T_{\Sigma} = \{f \in \al H(M,\bb C^n) : f \pitchfork \Sigma\}$ is open in $\al H(M,\bb C^n)$ with the weak topology .
\end{prop}

This case turned out to be too easy to handle.

Now, let $M$ be a Stein manifold and $N$ be any complex manifold. Let $\al H(M,N)$ denote the set of all holomorphic maps between $M$ and $N$. If we regard $M$ and $N$ as smooth manifolds and denote $C^{\infty}(M,N)$ the set of all smooth maps, then $C^{\infty}(M,N)$ can be given two topologies, the weak topology and the strong (Whitney) topology (when $M$ is compact, the two topologies are same), see Hirsch \cite{Hirsch}.  

Since $\al H(M,N) \subset C^{\infty}(M,N)$, we can give $\al H(M,N)$ the relative weak and strong topologies. However only the weak relative topology makes sense because the strong relative topology gives the discrete topology on $\al H(M,N)$. Thus we work with the weak topology on $\al H(M,N)$.

Many transversality theorems have been proved in this case (see Kaliman and Zaidenberg \cite{Kaliman} and Forstneri{\v{c}} \cite{Forstneric}). It is worth clarifying that in the transversality theorem (theorem 4.2) in \cite{Forstneric} we don't need the complex analytic subvarieties to be $a$-regular. Also, in lemma 4.4 in \cite{Forstneric} the complex analytic subvariety in the source manifold need not be $a$-regular and the one in the target manifold only has to be $a$-regular. 

Here we give a partial converse to the lemma 4.4 in  Forstneri{\v{c}} \cite{Forstneric} which says:

\begin{prop} Let $M$ (Stein) and $N$ be complex manifolds and let $B \subset N$ be a complex analytic subvariety of $N$. Let $\Sigma$ be an $a$-regular complex analytic stratification of $B$, then for any compact $K$ in $M$, the set $T_{\Sigma} = \{f \in \al H(M,N) : f \pitchfork_K \Sigma\}$ is open in the weak topology. 
\end{prop}

A partial converse of the above proposition reads:

\begin{prop} Let $M$ (Stein) and $N$ be complex manifolds and $\Sigma$ be a stratification of a complex analytic subvariety of $N$. Let $r = \min \{\dim S : S \,\,\text{is a stratum in}$ $\Sigma\}$. If $\dim M \geq \dim N - r$ and for some point $w \in M$, the set $T = \{f \in \al H(M,N) : f \pitchfork_w \Sigma\}$ is open in the weak topology, then $\Sigma$ is $a$-regular.
\end{prop}

\begin{pf} Without loss of generality, suppose that the strata are of dimension at least 1. Since $M$ is Stein, we can assume that $M$ is a complex submanifold of some $\bb C^p$ such that $w = {0} \in \bb C^p$, clearly $p \geq n - r$. 

We will again prove the proposition using contradiction. Suppose $\Sigma$ is not $a$-regular. Let $X$ and $Y$ be two distinct strata in $\Sigma$ such that $Y$ is not $a$-regular over $X$ at some $x \in X \cap \overline{Y}$. Then, there exists a sequence $\{y_i\}_{i=1}^{\infty}$ in $Y$ such that $\lim_{i \rightarrow \infty} y_i = x$, $\lim_{i \rightarrow \infty} T_{y_i} Y = \tau$ but $T_xX \not\subset \tau$. In fact, by the curve selection lemma, we can choose $y_i$ to lie on an analytic curve.

Let $v \in T_x X$ be such that $v \notin \tau$. Now, find a vector subspace $H$ of $\bb C^n$ of dimension $n - r$, as in proposition \ref{main},  such that
\begin{equation}
H + T_x X = T_x N 
 \end{equation}
and
\begin{equation}
 H + \tau \neq T_x N.
\end{equation}

Let $(\psi, V)$ be a coordinate chart around $x \in N$ such that $\psi(x) = \overrightarrow{0} \in \bb C^n$.  Then, $D_x\psi : T_xN \rightarrow \bb C^n$ is a linear isomorphism and under this isomorphism (2.1) and (2.2) become
\begin{equation}
D_x\psi(H) + D_x\psi(T_x X) = \bb C^n 
 \end{equation}
and
\begin{equation}
 D_x\psi(H) + D_x\psi(\tau) \neq \bb C^n.
\end{equation}

Let $m$ be the dimension of $M$ and choose a basis $\{{u_1},\ldots,{u_m},\ldots,{u_p}\}$ such that ${u_1},\ldots,{u_m}$ span the vector subspace $T_wM$ of $\bb C^p$.

Let $l$ be the dimension of $D_x\psi(H) + D_x\psi(\tau)$ ($n-r \leq l < n$). Now choose a basis $\{v_1,\ldots,v_l\}$ of  $D_x\psi(H) + D_x\psi(\tau)$ such that $\{v_1,\dots,v_{n-r}\}$ forms a basis of $D_x\psi(H)$and extend it to a basis $\{v_1,\ldots,v_l,v_{l+1},\ldots,v_{n-1},v'\}$ of $\bb C^n$ where $D_x\psi(v)=v'$.

Now, define a map $L : \bb C^p \rightarrow \bb C^n$ (this map is well defined because $p  > n-r$) by, 
$$L(a_1{u_1}+\cdots+a_p {u_p}) = a_1 {v_1} + a_2 {v_2} + \ldots + a_{n-r} {v_{n-r}}.$$
where ${v_1},\ldots,{v_{n-r}}$ is a basis of $D_x\psi(H)$.

Notice that, $L(T_wM) = D_x{\psi}(H)$ and $L$ is a holomorphic map. Let $g : \bb C^p \rightarrow N$ be defined by $g = \psi^{-1} \circ L$ then $g(w) = x$, $g$ is a holomorphic map and $D_wg = (D_x\psi)^{-1}L$ because $L$ is a linear map. Thus by (2.1), $g \pitchfork_w \Sigma$.

We will now construct a sequence of holomorphic maps $\{g^k\}$ between $\bb C^p$ and $N$ which converge to $g$ in the weak topology such that for sufficiently large $k$, $g_k \not\pitchfork_w \Sigma$.

As in proposition \ref{main}, first note that $y_k \in V$ for sufficiently large $k$. Thus, 
$$\lim_{k \rightarrow \infty} D_{y_k}\psi(T_{y_k} Y) =  D_x\psi(\tau).$$
Now, choose a basis $\{v_1^k,v_2^k,\ldots,v_l^k,v_{l+1}^k,\ldots,$ $v_{n-1}^k, v^k\}$ of $D_{y_k}\psi(T_{y_k}N)$ such that $D_{y_k}\psi(T_{y_k}Y)$ belongs to the span of $\{v^k_1,\ldots,v^k_l\}$ and
\begin{equation}
\lim_{k \rightarrow \infty} {v_i}^k = {v_i}\,\, \text{and}\,\,\lim_{k \rightarrow \infty} v^k = {v}
\end{equation}
Let $H^k$ be the subspace of $D_{y_k}\psi(T_{y_k}N)$ spanned by $v_1^k,\ldots,v_{n-r}^k$. Then, by (2.5) we have $\lim_{k\rightarrow \infty} H^k = D_x \psi(H)$ and we have 
\begin{equation}
H^k + D_{y_k}\psi(T_{y_k}Y) \neq \bb C^n.
\end{equation}
since the left hand side of (2.6) is spanned by a subset of $\{v_1,\ldots,v_l\}$ and $l < n$.
Now, define $L^k : \bb C^p \rightarrow \bb C^n$ by the formula
$$L^k(a_1{u_1}+\cdots+a_p {u_p}) = \psi(y_k) +a_1 {v_1^k} + a_2 {v_2^k} + \ldots + a_{n-r} {v^k_{n-r}}$$
and set $g^k = \psi^{-1} \circ L^k$. Clearly then, $g^k(w) = y_k$, by (2.5) the sequence of maps $g^k$ converges to the map $g$, $L^k(T_wM)=H^k$  and by (2.6), $g^k \not\pitchfork_{w} \Sigma$.

Now, let $g|_M = f$ and $g^k|_M = f^k$, then clearly $f$ and the $f^k$'s are holomorphic maps between $M$ and $N$, the sequence $f^k$ converges to $f$ in the weak topology and $f \pitchfork_w \Sigma$ but for large enough $k$, $f^k \not\pitchfork_w \Sigma$, which is a contradiction to the hypothesis that the set $T$ is open in the weak topology.
\end{pf}

\section*{Acknowledgements}
I would like to thank my thesis advisor Professor David Trotman for helpful discussions, suggestions, comments, encouragement and careful reading of the paper. I would also like to thank Laurent Battisti for discussions clarifying some basic ideas of complex analysis.

\nocite{Ricketts}

\bibliographystyle{amsplain}
\bibliography{bibliography}

\end{document}